\documentclass[12pt]{amsart}
\usepackage{amscd}
\theoremstyle{definition}

\def \ph{\varphi}
\def \ra{\rightarrow}

\def \hom{\mbox{\rm Hom}}

\def \tns{\otimes}

\def \mcom{,\cdots,}
\def \k{\mbox{$\mathfrak K$}}

\def \Z{\mbox{$\mathbb Z$}}

\def\br#1#2{\lbrack#1,#2\rbrack}

\def\zt{\mbox{$\Z_2$}}

\def\sh{\operatorname{Sh}}
\def\inv{^{-1}}
\def\d{d}

\def\A{\mbox{$\mathcal A$}}
\def\L{L}
\def\m{\mbox{$\mathfrak m$}}

\def\dl{\delta}

\def\coder{\operatorname{Coder}}

\def\linf{\mbox{$L_\infty$}}

\def\s#1{(-1)^{#1}}
\def\tl{\tilde\lambda}
\def\HL{H(\L)}
\def\ZL{Z(\L)}

\def\SW{S(W)}
\def\H{\mathcal H}
\def\th#1#2{\theta^{#1}_{#2}}
\def\pha#1#2{\ph^{#1}_{#2}}
\def\ps#1{\psi_{#1}}
\def\dl#1#2{\delta^{#1}_{#2}}
\def\psa#1#2{\psi^{#1}_{#2}}
\def\ts#1#2{t^{#1}_{#2}}
\def\tth#1#2{{\tilde{\theta}}^{#1}_{#2}}
\def\tph#1#2{{\tilde{\ph}}^{#1}_{#2}}
\def\tps#1#2{{\tilde{\psi}}^{#1}_#2}
\author{Alice Fialowski}
\address{E\"otv\"os Lor\'and University\\
Budapest, Hungary}
\email{fialowsk@cs.elte.hu}
\author{Michael Penkava}
\address{University of Wisconsin\\
Eau Claire, WI 54702-4004}
\email{penkavmr@uwec.edu}
\thanks{The research of the authors was partially supported by
grants OTKA T030823, T29525, FKFP 0170/1999 and by
grants from the University of Wisconsin-Eau Claire}
\title[Infinity Algebras and Deformations]{Examples of Infinity and Lie Algebras and Their Versal Deformations}
\begin{document}
\setlength{\multlinegap}{0pt}
%\nocite{ps2,pen1,ls,pen2,pen3,kon,
%fm,mar,mar2,ksv,sta1,sta2,sta3,getz,getz2,ge_ka1,ge_ka2,
%gers,lod,hoch,fi,fi2,ff2,ff,ff3}
%\input abstract.tex
\begin{abstract}
This article explores some simple examples of $L_\infty$ algebras and
the construction of miniversal deformations of these
structures. Among other things, it is shown that
there are two families of nonequivalent $L_\infty$ structures
on a $1|1$ dimensional vector space, two of which are Lie algebra
structures. The main purpose of this work is to provide a simple
effective procedure for constructing miniversal deformations,
using the examples to illustrate the general technique.  The
same method can be applied directly to construct versal deformations of
Lie algebras.
\end{abstract}
\date{January 22, 2001}
\maketitle
%\table
%\input intro.tex
\section{Introduction}

In this article we shall explore some interesting simple examples of
\linf\ structures on low dimensional \zt-graded vector spaces, and
construct versal deformations of these structures. A detailed
explanation of the definitions of infinity algebras and their cohomology
theory can be obtained in \cite{pen1}. Infinity algebras
were first described and applied by James Stasheff (see
\cite{sta1,sta2,ls,lm}), and have appeared recently in both mathematics
and mathematical physics (see
\cite{aksz,bfls,kon,mar2,sta3,lp,ps2,pen3}). The essential details of the
theory of versal deformations of infinity algebras are outlined in
\cite{fiapen1}, and are a generalization of the results in \cite{ff2}
on versal deformations of Lie algebras (see also \cite{fi,fi2,ff3,fp,ss}).
In this article, our aim is
to show how to construct a miniversal deformation through some examples.
For the proof of the miniversality of the deformations we will construct
here, we refer to \cite{fiapen1}. For sake of completeness, we will
present here the basic definitions of \linf\ algebras and formal
deformations.

If $W$ is a \zt-graded vector space, then $\SW$ denotes the symmetric
coalgebra of $W$. If we let $T(W)$ be the reduced tensor algebra
$T(W)=\bigoplus_{n=1}^\infty W^{\tns n}$, then the reduced symmetric algebra
$S(W)$ is the quotient of the tensor algebra by the graded ideal
generated by $u\tns v-\s{uv}v\tns u$ for elements $u$, $v\in W$.
The symmetric algebra has a natural coalgebra structure, given by
\begin{equation*}
\Delta(w_1\dots w_n)=
\sum_{k=1}^{n=1}\sum_{\sigma\in\sh(k,n-k)}\epsilon(\sigma)
w_{\sigma(1)}\dots w_{\sigma(k)}\tns w_{\sigma(k+1)}\dots w_{\sigma(n)},
\end{equation*}
where we denote the product in $S(W)$ by juxtaposition,
$\sh(k,n-k)$ is the set of \emph{unshuffles} of type $(k,n-k)$,
and $\epsilon(\sigma)$ is a sign determined by $\sigma$ (and $w_1\dots w_n$)
given by
$$
w_{\sigma(1)}\dots w_{\sigma(n)}=\epsilon(\sigma)w_1\dots w_n.
$$
A coderivation on $S(W)$ is a map $\delta:S(W)\ra S(W)$ satisfying
$$\Delta\circ\delta=(\delta\tns I+I\tns\delta)\circ\Delta.$$
Let us suppose that the even part of $W$ has basis
$e_1\dots e_m$, and the odd part has basis $f_1\dots f_n$, so that $W$
is an $m|n$ dimensional space. Then a basis of $\SW$ is given by all
vectors of the form $e_1^{k_1}\dots e_m^{k_m}f_1^{l_1}\dots f_n^{l_n}$,
where $k_i$ is any nonnegative integer, and $l_i\in\zt$. An \linf\
structure on $W$ is simply an odd codifferential on $\SW$, that is to
say, an odd coderivation whose square is zero. The space $\coder(W)$
can be naturally identified with $\hom(S(W),W)$, and the Lie
superalgebra structure on $\coder(W)$ determines a Lie bracket on
$\hom(S(W),W)$ as follows. Denote $L_m=\hom(S^m(W),W)$ so that
$L=\hom(S(W),W)$ is the direct product of the spaces $L_i$.
If $\alpha\in L_m$ and $\beta\in L_n$, then $\br\alpha\beta$ is the
element in $L_{m+n-1}$ determined by
\begin{multline}\label{braform}
\br\alpha\beta
(w_1\dots w_{m+n-1})=\\
\sum_{\sigma\in\sh(n,m-1)}\epsilon(\sigma)\alpha(\beta(w_{\sigma(1)}\dots
w_{\sigma(n)})w_{\sigma(n+1)}\dots w_{\sigma(m+n-1)})\\
-\s{\alpha\beta}\
\sum_{\sigma\in\sh(m,n-1)}\epsilon(\sigma)
\beta(\alpha(w_{\sigma(1)}\dots
w_{\sigma(m)})w_{\sigma(m+1)}\dots w_{\sigma(m+n-1)}).
\end{multline}
Another way to express this bracket is in the form
\begin{equation*}
\br\alpha\beta=\alpha\tilde\beta-\s{\alpha\beta}\beta\tilde\alpha,
\end{equation*}
 where for $\ph\in\hom(S^k(W),W)$, $\tilde\ph$ is the
associated coderivation, given by
\begin{equation*}
\tilde\ph(w_1\dots w_n)=
\sum_{\sigma\in\sh(k,n-k)}\epsilon(\sigma)
\ph(w_{\sigma(1)}\dots w_{\sigma(k)})w_{\sigma(k+1)}\dots w_{\sigma(n)}.
\end{equation*}
In this article, formula (\ref{braform}) will be used extensively
to make explicit computations of brackets.  In our computations,
we shall also use the fact that there are $n\choose k$  permutations in
$\sh(k,n-k)$.

If $W$ is completely odd, and $d\in\L_2$, then $d$ determines an ordinary
Lie algebra on $W$, or rather on its parity reversion. The
symmetric algebra on $W$ looks like the exterior algebra on $W$ if
we forget the grading.  If we define $[a,b]=d(ab)$ for $a,b\in W$, then
the bracket is antisymmetric because $ba=-ab$, and moreover
\begin{align*}
0=\br dd(abc)=&\
\frac12\sum_{\sigma\in\sh(2,1)}\epsilon(\sigma)
d(d(\sigma(a)\sigma(b))\sigma(c))\\=&
d((d(ab)c)+d(d(bc)a)-d(d(ac)b)\\
=&[[a,b],c]+[[b,c],a]-[[a,c],b],
\end{align*}
which is the Jacobi identity.  When $d\in\L_2$ and $W$ has a true grading,
then the same principle holds, except that one has to take into
account a sign arising from the map $S^2(W)\ra\bigwedge^2(V)$, where
$V$ is the parity reversion of $W$.  Thus \zt-graded Lie algebras
are also examples of \linf\ algebras.  In addition, differential
graded Lie algebras are examples of \linf\ algebras.  In all these
cases, the method of construction  of miniversal deformations we
describe here applies.  One simply considers only terms that
come from $\L_2$ in the codifferentials.

Suppose that $\tilde g:S(W)\ra S(W')$ is a coalgebra morphism,
that is a map satisfying
$$
\Delta'\circ\tilde g=(\tilde g\tns \tilde g)\circ\Delta.
$$
If $d$ and $d'$ are \linf\ algebra structures on $W$ and $W'$,
resp., then $\tilde g$ is a homomorphism between these structures if
$\tilde g\circ d=\d'\circ\tilde g$. Two \linf\ structures $d$ and $d'$ on $W$
are equivalent when there is a coalgebra automorphism $\tilde g$ of $S(W)$
such that $d'={\tilde g}\inv\circ d\circ {\tilde g}$.  Furthermore,
if $d=d'$, then $\tilde g$ is said to be an automorphism of the \linf\
algebra.
\section{Deformations with a local base}
\def\da{d_{\A}}
We shall always assume that the ground field $\k$ of the
graded space $W$ has characteristic zero. A local base
is a \zt-graded commutative local algebra $\A$ equipped with
a fixed augmentation $\epsilon:\A\ra\k$, with augmentation ideal $\m$.
If $d$ is a codifferential on $S(W)$, then a deformation of $d$ with
base $\A$ is a codifferential $\da$ on $S(W)\tns\A$, which projects
to $d$ under the morphism of coalgebras $I\tns\epsilon:S(W)\tns\A\ra S(W)$.
When $\A$ is a formal algebra (a complete local algebra), then a
deformation with base $\A$ will be called a formal deformation, and
a deformation with base $\A/\m^{n+1}$ will be called an $n$-th order
deformation.  When $n=1$, the deformation is called infinitesimal.

In \cite{fiapen1}, it was shown that under some mild restrictions on the
cohomology of $L$, there is a miniversal deformation,
that is, there is a certain local base $\A$ and deformation $\da$,
such that if $\da'$ is any deformation with base $\A'$, then there
is a morphism of algebras $\lambda:\A\ra \A'$, such that
the induced codifferential $\lambda_*(\da)$ on $S(W)\tns\A'$
is equivalent to $\da'$. Note that any coderivation of $S(W)\tns \A$
lies in $L\tns \A$, so that $\lambda_*$ is simply the natural
map $L\tns \A\ra L\tns\A'$. Moreover, the notion of equivalence
of deformations with base $\A$ is given by the following.
If $d$ is a codifferential on $S(W)$ and $\tilde g$ is a coalgebra
automorphism of $S(W)$, then define
${\tilde g}^*(d)={\tilde g}\inv\circ d\circ \tilde g$.
If $\da^1$ and $\da^2$ are two deformations of $d$ with base $\A$,
then they are said to be equivalent if there is a coalgebra
automorphism $\tilde g$ of $S(W)\tns \A$ satisfying
$(I\tns\epsilon)\circ \tilde g=I\tns\epsilon$ such that
${\tilde g}^*(\da^1)=\da^2$.
The condition on $\tilde g$ means simply that $\tilde g(w)
=w\pmod{\m}$, so that $\tilde g$ does not act on the underlying coalgebra
$S(W)$ (to avoid mixing the notion of automorphisms of the \linf\
algebra with that of the equivalences of deformations of the algebra).

The objective in this article is to
compute versal deformations of some simple examples.
To see how to do this in practice, we will give a concrete description
of the construction of a miniversal deformation of a codifferential $d$.
Since $d$ is a codifferential, the map $D:L\ra\L$ given by $D(\ph)=[d,\ph]$
satisfies $D^2=0$. The homology $H(L)$ of this map is called the cohomology
of the \linf\ algebra.
Let us assume that
the space of cocycles $\ZL$ has an increasing
basis of the form $\langle \delta_i,\beta_j\rangle$,
where the subspace spanned by the $\delta_i$ projects isomorphically onto
the cohomology of $d$, and $\beta_j$ is a basis of the coboundaries.  By
basis, we mean (in the sense of \cite{fiapen1})
that any element of $\ZL$ has a unique expression as a power series in
the elements of the basis, and increasing means that the orders of
the elements is nondecreasing.  Here, order refers to the order of the
element in $L$.  In the general case, an element in $L$ is represented
as a power series, since $L$ is a direct product of the spaces $L_n$. The
order represents the least $n$ such that the projection of the element
onto $L_n$ is nonzero.

In order for there to be a good
theory of deformations, the cohomology should be of finite type,
which means that there are only a finite number of elements in
the basis $\delta_i$ of order less than a given integer.  In our examples,
the condition will be satisfied, but even in the simplest cases,
as we shall see in the case of a one dimensional space, the cohomology
can be infinite dimensional.

To construct the universal infinitesimal deformation of a codifferential
$d$, one formally considers the infinitesimal base $\A=\k\oplus \H$
where $\H=(\Pi(\HL))^*$ is the dual of the parity reversion of the
cohomology of $L$, and the multiplication in $\H$ is trivial.
Then the base of the miniversal deformation is $\A=\k[[\H]]/R$,
where $R$ represents a set of polynomial relations involving the
generators of the space $\H$.

If we let $\delta_i$ be a basis of a subspace of cocycles which projects
isomorphically to the cohomology, then choosing an appropriate dual
basis $u^i$ of $\H$, the universal infinitesimal deformation $d_1$
is given by
\begin{equation}
d_1=d+\delta_iu^i.
\end{equation}
The parity of the element $u^i$ is opposite to the parity of $\delta_i$,
in order to preserve the oddness of the codifferential.
Let $\alpha_1=\delta_iu^i$. Then  $\alpha_1$ is a cocycle.
In order to
extend the bracket to a second order bracket, one simply computes
\begin{equation}
\br{d_1}{d_1}=\br dd +2\br d{\alpha_1}+\br{\alpha_1}{\alpha_1}
=\br{\alpha_1}{\alpha_1}
\end{equation}
which is a cocycle, because the bracket of two cocycles is always a cocycle.
This cocycle can be represented as a linear combination of the $\delta_i$
plus a coboundary term.
Let us express this in the form
\begin{equation}
\br{d_1}{d_1}=-\frac12D(\alpha_2)+\delta_iR^i_2,
\end{equation}
where $R^i_2$ is a sum of products of the parameters of the form $u_ku_l$,
so lies in $\m^2$,
and $\alpha_2\in L\tns \m^2$.
We would like $\br{d_1}{d_1}$ to be a coboundary mod $\m^3$, which means
that we must have $R^i_2=0$ mod $\m^3$.

The second order deformation
is given by $d_2=d_1+\alpha_2$.
The bracket of this second order codifferential with itself is
\begin{align*}
\br{d_2}{d_2}&=
\br{d_1}{d_1}+2\br{d_1}{\alpha_2}+\br{\alpha_2}{\alpha_2}\\
&=
\br{d_1}{d_1}+2D(\alpha_2)+2\br{\alpha_1}{\alpha_2}+\br{\alpha_2}{\alpha_2}\\
&=\delta_iR^i_2+ 2\br{\alpha_1}{\alpha_2}+\br{\alpha_2}{\alpha_2}
\end{align*}
Note that the first and third terms in this bracket are zero mod $\m^3$.
We want this bracket to be a coboundary mod $\m^3$. Let us first show
that at least it is a cocycle mod $\m^3$.
Now
\begin{align*}
D(\br{d_2}{d_2})&=
2\br{\alpha_1}{D(\alpha_2)}+2\br{D(\alpha_2)}{\alpha_2}\\
&=2\br{\alpha_1}{-2\br{\alpha_1}{\alpha_1}+2\delta_iR^i_2}
+2\br{D(\alpha_2)}{\alpha_2}\\
&=4\br{\alpha_1}{\delta_iR^i_2}
+2\br{D(\alpha_2)}{\alpha_2},
\end{align*}
which is equal to zero mod $\m^4$.
In general, let us suppose that we have shown that $\br{d_n}{d_n}=0$
mod $\m^{n+1}$. and that $D(\br{d_n}{d_n})=0$ mod $\m^{n+2}$.
Then we set
\begin{equation*}
\br{d_n}{d_n}=-\frac12D(\alpha_{n+1})+\delta_iR^i_{n+1}
\end{equation*}
where $\alpha_{n+1}\in L\tns\m^{n+1}$, and $R^i_{n+1}=0$ mod $\m^{n+2}$.
Set $d_{n+1}=d_n+\alpha_{n+1}$.
Then
\begin{align*}
\br{d_{n+1}}{d_{n+1}}&=
\delta_iR^i_{n+1}+2\br{d_n-d}{\alpha_{n+1}}+\br{\alpha_{n+2}}{\alpha_{n+2}}.
\end{align*}
Note that each term appearing in the bracket is equal to zero mod $\m^{n+2}$.
Next, we compute
\begin{align*}
D\br{d_{n+1}}{d_{n+1}}=&
2\br{D(d_n-d)}{\alpha_{n+1}}-
2\br{d_n-d}{D\alpha_{n+1}}+2\br{\alpha_{n+1}}{\alpha_{n+1}}\\
=&
2\br{D(d_n-d-\alpha_1)}{\alpha_{n+1}}-
2\br{d_n}{D\alpha_{n+1}}+2\br{\alpha_{n+1}}{\alpha_{n+1}}\\
=&
2\br{D(d_n-d-\alpha_1)}{\alpha_{n+1}}-
2\br{d_n}{-2\br{d_n}{d_n}+2\delta_iR^i_{n+1}}\\&+
2\br{\alpha_{n+1}}{\alpha_{n+1}}\\
=&
2\br{D(d_n-d-\alpha_1)}{\alpha_{n+1}}-
2\br{d_n,2\delta_iR^i_{n+1}}+
2\br{\alpha_{n+1}}{\alpha_{n+1}}\\
=&
2\br{D(d_n-d-\alpha_1)}{\alpha_{n+1}}-
2\br{d_n-d}{2\delta_iR^i_{n+1}}\\&+
2\br{\alpha_{n+1}}{\alpha_{n+1}}.
\end{align*}
Each of the three remaining terms is equal to zero mod $\m^{n+3}$, so
one can continue the process by adding an $\alpha_{n+2}$ term. Thus
we see that it is possible to construct a formal deformation
$d_\infty=\sum_{i=1}^\infty{\alpha_i}$, which will satisfy
$\br{d_\infty}{\d_\infty}=0$. But what about the relations $R^i_{n}$?
It is easy to see that $R^i_{n+1}=R^i_n$ mod $\m^{n+1}$. We obtain
a series of relations $R^i_\infty$, which must hold on the base in
order for the required bracket relation to hold.

When it happens that  $\br{d_n}{d_n}$ contains no
coboundary terms for some $n$, then the $n$-th order deformation is miniversal,
and the $n$-th order relations hold to all orders.  Thus, in good
cases, one may hope to calculate the miniversal deformation exactly.
This is what happens in the examples we shall study.

We shall only consider the simplest cases of \linf\ algebras. It
is interesting to note that even for a 2 dimensional vector space,
it can happen that there are infinitely many nonequivalent \linf\
structures.
The main purpose of this article
is to demonstrate how to compute versal deformations in practice.

Note
that in order to calculate a versal deformation, we need much
more information than the dimensions of the cohomology groups.
Thus, calculations of cohomology of Lie algebras which
often appear in the literature are not concrete enough to do
versal deformations, because the bracket structure on the cohomology
is essential for the calculations, because only the dimensions of the
cohomology groups are computed.  In fact in the examples presented
in this paper, we had to also compute
brackets of cochains that are not cocycles.
\section{Deformations of a $0|1$ dimensional space}
Suppose that $W=\langle f\rangle$ is a 1 dimensional odd vector space.
This situation corresponds by parity reversion to the case of an
ordinary even one dimensional space, on which there is obviously
only the trivial Lie algebra structure.  Indeed, there is only
one \linf\ algebra structure $d=0$ as well, so it does not seem
reasonable to expect any deformations of this trivial structure.
Nevertheless, this intuition is slightly wrong, because there is
a infinitesimal differential on the space $W$, and a differential
graded space is an \linf\ algebra. To see what is going on, note
that $L_1=\langle \ph\rangle$, where $\ph(e)=e$ is the identity
map on $W$, which is an even 1-cochain. The map $\ph$ does not determine
an \linf\ structure, because it is even, but when multiplied by
the odd parameter $\theta$, it becomes an odd codifferential on
$S(W)=W$. Thus $d_1=\ph\theta$ is an infinitesimal deformation of
the trivial deformation $d=0$.  Since $\ph$ is a cocycle, but not
a coboundary, this deformation is also a miniversal deformation.
To determine any relations on the base $\k[[\theta]]=\k\oplus\k\theta$,
one should compute $\br{d_1}{d_1}$, but this is automatically
zero since $\theta^2=0$. Thus, there are no relations other than
the relation $\theta^2=0$, which is always true for an odd parameter,
and the base of the miniversal deformation is just $A=\k[\theta]$.
\section{Deformations of a $1|0$ dimensional space}
Let $W=\langle e\rangle$ be a one dimensional even vector space.
Then $S^k(W)=\langle e^k\rangle$, and $L_k=\langle \ph_k\rangle$,
where $\ph_k(e^l)=k!\dl lk e$. Now $\br{\ph_k}{\ph_l}$ is the element
of $L_{k+l-1}$, which we compute in detail, as an illustration of how
to perform the computation of brackets. Let $n=k+l-1$. Then
\begin{multline*}
\br{\ph_k}{\ph_l}(e^n)=\\
\sum_{\sigma\in\sh(l,n-l)}k!l!
\ph_k(\ph_l(e^l),e^{k-1})
-\s{\ph_k\ph_l}
\sum_{\sigma\in\sh(k,l-k)}l!k!\ph_k(\ph_k(e^k),e^{l-1})\\
=\left(k!l!{n\choose l}-k!l!{n\choose k}\right)e
=(k-l)n!e=(k-l)\ph_n(e^n).
\end{multline*}
Since all the maps $\ph_k$ are even, the only codifferential on
$W$ is $d=0$. This time, the infinitesimal deformation
$d_1=\ph_k\theta^k$ has an infinite number of terms, and they
are not coboundaries. There are some nontrivial
relations on the base, arising from the self bracket of $d_1$.
We have
\begin{equation*}
\frac12\br{d_1}{d_1}=
\sum_{n=1}^\infty
\sum_{k+l=n+1}\br{\ph_k\theta^k}{\ph_l\theta^l}
=
\sum_{n=1}^\infty
\sum_{k+l=n+1}(k-l)\theta^k\theta^l\ph_n,
\end{equation*}
from which it follows that we obtain an infinite set of relations of the
form
$\sum_{k+l=n+1}(k-l)\theta^k\theta^l=0$, for $n=1\dots$.
Let us examine the first few terms.  For $n=1$, the relation is zero.
For $n=2$, we obtain
$\theta^2\theta^1-\theta^1\theta^2=0$, i.e. $\theta^1\theta^2=0$.
For $n=3$, the relation reduces to $\theta^1\theta^3=0$. For $n=4$,
we obtain $6\theta^1\theta^4+10\theta^2\theta^3=0$. If we let
$R=\langle
\sum_{k+l=n+1}(k-l)\theta^k\theta^l=0
\rangle$
be the ideal generated by these relations in $\k[[\theta^1\dots]]$,
then the base of the miniversal deformation is just
$\A=\k[[\theta^1\dots]]/R$. The miniversal deformation coincides
with the universal infinitesimal deformation, except that the base
$\A$ of the miniversal deformation is not infinitesimal, but is
a somewhat complicated infinite dimensional graded commutative algebra.
\section{Deformations of a $2|0$ dimensional space}
Let $W=\langle e_1, e_2\rangle$ be a completely even dimensional
space. Then $S^k(W)=\langle e_1^k, e_1^{k-1}e_2\mcom e_2^k\rangle$
has dimension $k+1$, so that $L_k$ has dimension $2(k+1)$. It will
simplify the notation if we introduce multi-indices of nonnegative
integers $I=(i_1,i_2)$. Define
$|I|=i_1+i_2$, $I!=i_1!i_2!$, $\dl IJ=\dl {i_1}{j_1}\dl {i_2}{j_2}$ and
$e^I=e_1^{i_1}e_2^{i_2}$.
Then $L_n=\langle \ph_{I,i}:|I|=n, i\in\{1,2\}\rangle$, where
\begin{equation*}
\ph_{I,i}(e^{J})=I!\dl JI e_i.
\end{equation*}
Note that because  $W$ is totally even, the maps $\ph_{I,i}$ are all
even, so that $d=0$ is the only \linf\ structure on $W$ and
$$d_1=\ph_{I,i}\theta^{I,i}$$ is the universal infinitesimal deformation
as well as a miniversal deformation, but in the latter case, we
need to compute the relations on the base. To that end, let us compute
$\br{\ph_{I,k}}{\ph_{J,l}}$, which lies in $L_{|I|+|J|-1}$.
Define $I+J=(i_1+j_1,i_2+j_2)$ and
$$
I-k=
\begin{cases}
(i_1-1,i_2)& \hbox{if $k=1$ and $i_1>0$}\\
(i_1,i_2-1)& \hbox{if $k=2$ and $i_2>0$}\\
\end{cases}
$$
Then
\begin{equation*}
\br{\ph_{I,k}}{\ph_{J,l}}=
i_l\ph_{J+(I-l),k}-
j_k\ph_{I+(J-k),l},
\end{equation*}
where terms such that $I-l$ or $J-k$ are not defined are omitted.
 From this we determine that the relations are given by
\begin{equation*}
\sum_{J+(I-j)=M}i_j\theta^{I,m}\theta^{J,j}=0,
\end{equation*}
for each multi-index $M$ and $m\in\{1,2\}$.
This follows from the requirement that $\br{d_1}{d_1}=0$.
Computing this bracket, we obtain,
\begin{multline*}
\frac12\br{d_1}{d_1}=
(i_l\ph_{J+(I-l),k}-
j_k\ph_{I+(J-k),l})\theta^{I,k}\theta^{J,l}
=\\
\sum_{J+(I-l)=M\atop k=m}i_l\ph_{M,m}\theta^{I,m}\theta^{J,j}
-
\sum_{I+(J-l)=M\atop l=m}j_k\ph_{M,m}\theta^{I,k}\theta^{J,m}.
\end{multline*}
The second term is the same as the first after interchanging the
dummy indices and using the fact that the $\theta$'s are odd.
\section{Deformations of a $0|2$ dimensional space}

Let $W=\langle f_1,f_2\rangle$ be a completely odd two dimensional
space. Then $S^2(W)=\langle f_1f_2\rangle$, and $S(W)=W+S^2(W)$.
We have $L_1=\langle \pha ij\rangle$, where $\pha ij(f_k)=\delta^i_kf_j$,
and $L_2=\langle \ps i\rangle$, where $\ps i(f_1f_2)=f_i$.
Note that $L_1$ consists of even elements, while $L_2$ is odd. This
time, in addition to the trivial \linf\ structure $d=0$, any element
of $L_2$ also determines an \linf\ algebra structure.  All of these
structures are actually Lie algebra structures, with $d=0$ giving the
abelian Lie algebra structure, while the nonzero elements of $L_2$
give rise to equivalent structures, so we need only consider the case
$d=\ps 1$.

As an aid to the construction, we first compute the brackets of all the
basis elements. First,
\begin{equation*}
\br{\pha ij}{\pha kl}(f_m)=
\pha ij(\pha kl(f_m))-\pha kl(\pha ij(f_m))=
\delta^i_l\delta^k_mf_j-\delta^k_j\delta^i_m f_l,
\end{equation*}
while
\begin{multline*}
[\pha ij,\ps k](f_1f_2)=
\pha ij(\ps k(f_1f_2))-
\ps k(\pha ij(f_1)f_2)+\ps k(\pha ij(f_2)f_1)=\\
\pha ij(f_k)-\ps k(\delta^i_1f_jf_2)+\ps k(\delta^i_2f_jf_1)=
\delta^i_kf_j-\delta^i_1\delta^1_jf_k-\delta^i_2\delta^2_jf_k,
\end{multline*}
so that the brackets are given by
\begin{align*}
[\pha ij,\pha kl]&=\delta^i_l\pha kj-\delta^k_j\pha il\\
[\pha ij,\ps k]&=\delta^i_k\ps j-\delta^i_j\ps k\\
[\ps i,\ps j]&=0.
\end{align*}
We give a more detailed table of the brackets below, which is useful in determining
the cocycles and coboundaries.
\begin{equation*}
\begin{array}{llll}
\br{\pha 11}{\pha11}=0\\
\br{\pha 11}{\pha12}=-\pha12\\
\br{\pha 11}{\pha21}=\pha21&\br{\pha12}{\pha21}=\pha22-\pha11\\
\br{\pha 11}{\pha22}=0&\br{\pha12}{\pha22}=-
\pha12&\br{\pha21}{\pha22}=\pha21\\
\br{\pha 11}{\ps 1}=0&\br{\pha12}{\ps 1}=\ps 2&\br{\pha21}{\ps 1}=0&\br{\pha22}{\ps 1}=-
\ps 1\\
\br{\pha 11}{\ps 2}=-\ps 2&
\br{\pha 12}{\ps 2}=0&
\br{\pha21}{\ps 2}=\ps 1&
\br{\pha22}{\ps 2}=0
\end{array}
\end{equation*}
\subsection{Case 1: $d=0$}

For the trivial codifferential $d=0$, every cochain is a cocycle, and
therefore, the universal infinitesimal deformation of $d$ is
\begin{equation*}
d_1=\pha ij\th ji+\ps kt^k,
\end{equation*}
where as usual $\th ji$ are odd parameters and $t_k$ are
even parameters. Computing the bracket of this derivation with itself,
one obtains
\begin{equation*}
\frac12\br{d_1}{d_1}=
-(\delta^i_k\ps j+\delta^i_j\ps k)\th ji t^k
+(\delta^i_n\ph^m_j-\delta^n_j\ph^i_m)\th ji\th mn.
\end{equation*}
Since this bracket must vanish for $\d_1$ to be a codifferential,
we obtain some relations on the parameters. These are more easily
seen by examining the detailed table, and looking for which terms
contribute to which output cochains. For example, the cochain $\ps 1$
only arises from the brackets $\br{\ph^2_2}{\ps 1}$ and
$\br{\ph^2_1}{\ps 2}$, and therefore, $\th 22 t^1-\th 21t^2=0$.
The complete set of relations is given by
\begin{equation*}
R=\{
\th11\th211+\th21\th22,
\th11\th12+\th12\th22,
\th12\th21,
\th11t^2-\th21t^1,
\th22t^1-\th21t^2
\},
\end{equation*}
and thus the base of the miniversal deformation is
$\A=\k[[\theta^i_j,t_k]]/R$.
\subsection{Case $d=\ps 1$}

This second case is the first nontrivial codifferential that we have
encountered. From the table of brackets, we see that
$\{\ph^1_1,\ph^2_1,\ps 1,\ps 2\}$ is a basis of the space of cocycles,
with $\{\ps 1,\ps 2\}$ being a basis of the coboundaries.  Thus
$\{\ph^1_1,\ph^2_1\}$ projects to a basis of the cohomology of $d$,
which is two dimensional.
The universal infinitesimal deformation of $d$ is given by
$d_1=\ps 1+\ph^1_1\th11+\pha21\th12$. Since
$\br{d_1}{d_1}=2\pha21\th11\th12$, we obtain only one
relation, $R=\{\th11\th12\}$, so that the infinitesimal
deformation is miniversal, and
$\A=\k[\th11,\th12]/R=\k+\k\th11+\k\th12$ is the
base of the versal deformation.
This completes the description of the miniversal deformation;
its formula coincides with that of the
universal infinitesimal deformation. We would like  to
explore this example a bit further.

We will consider a different
infinitesimal deformation, which we will extend to a formal deformation.
Since the infinitesimal deformation we have already given is universal,
we start by adding some coboundary terms to our original deformation.
We begin with
$$d_1'=\ps 1+\pha11\th11+\pha21\th12+\ps1t^1+\ps2t^2,$$
and proceed to construct a formal deformation out of this infinitesimal
deformation of $d$.  Since we added merely a coboundary term, this deformation is
equivalent in some sense to the original one, but note that the
base of this infinitesimal deformation is different, since there are more parameters.

We first carry out the construction of a formal deformation from the
above infinitesimal one, and then show how to construct the homomorphism
from the base of the miniversal deformation to the base of the formal
deformation we are about to construct. Then we will explicitly determine
the equivalence between the push out of $d$ given by the homomorphism
of the bases and the formal deformation $d_\infty'$ which is obtained by
extending the infinitesimal $d_1'$ to a formal one.  The authors found
this construction intriguing, because the obvious homomorphism between
the bases does not turn out to be the correct one. In addition, this example
illustrates how to extend a deformation to higher order by adding
higher order terms to eliminate coboundary terms arising in the computation of
the self-bracket of the codifferential. Up to now, all of the terms
arising in these computations have been coboundary free, so that they
generate only relations, which means that the associated formal
deformations are really infinitesimal ones and thus don't provide a
good illustration of the theory of formal deformations.

Let us first compute the bracket $\br{d_1'}{d_1'}$ of the infinitesimal
deformation.
We have
\begin{align*}
\frac12\br{d_1'}{d_1'}&=
\br{\pha11}{\pha21}\th11\th12-\br{\pha11}{\ps2}\th11t^2-
\br{\pha21}{\ps2}\th12t^2\\
&=\pha21\th11\th12+\ps2\th11t^2-\ps1\th12t^2\\
&=\pha21\th11\th12-D(\pha12)\th11t^2-D(\pha22)\th12t^2.
\end{align*}
The second order deformation is given by extending $d_1'$ to
\begin{equation*}
d_2'=
\ps1+\pha11\th11+\pha21\th12+\ps1t^1+\ps2t^2+\pha12\th11t^2+\pha22\th12t^2.
\end{equation*}
There is one necessary relation on the parameters: $\th11\th12=0$.  Note
that formally speaking, this relation should be interpreted as being
only up to order 2, so that in the formal deformation, the relation that
corresponds to this one could, a priori, pick up some higher order
terms.

One could continue in the same manner to compute a third order
deformation, but note that if the expression above is bracketed with
itself, since $\ps2$ kills $\pha12$ and $\pha22$, and $\ps1$ takes
$\pha12$ to $\ps2$ and $\pha22$ to $\ps1$, we will end up only adding more
terms involving $\pha12$ and $\pha22$.  Examining the terms in the bracket which
involve $\pha11$ and $\pha21$ one determines that the same relation
$\th11\th12=0$ holds (up to third order) and no additional relations
are necessary.
Thus it is reasonable to expect that one can write
the formal deformation in the form
\begin{equation*}
\d_\infty'=
\ps1(1+t^1)+
\pha11\th11+\pha21\th12+\ps2t^2+\pha12\th11 t^2c_1+\pha22\th12t^2c_2,
\end{equation*}
and solve for the constants $c_1,c_2$  and the relations on the parameters
necessary so that
$\br{\d_\infty'}{\d_\infty'}=0$.
The solution is given by $c_1=c_2=\frac1{1+t^1}$,
which should be considered as a formal power series in $t^1$. The base
of this formal deformation is
$\A'=\k[[\th12,\th11,t^1,t^2]]/(\th11\th12)$.

Since the deformation $d_1$ is miniversal, there must be some morphism
$\lambda:\A\ra \A'$, such that $\lambda_*(d_1)\cong d_\infty'$. This means that
there is some automorphism $\tilde g$ of $S(W)$ such that
$ d_\infty'={\tilde g}\inv\lambda_*(d_1){\tilde g}$. Suppose that $\tilde
g$ is determined by an isomorphism $g:W\tns\A'\ra W\tns\A'$ of the form
$g=I+u_1\pha12 +u_2\pha22$, for some constants $u_1$ and $u_2$. First
we show that this is indeed an invertible map, by constructing its
inverse $g\inv=I+v_1\pha12 +v_2\pha22$. It is easily checked  that
$v_1=\frac{-u_1}{1+u_2},v_2=\frac{-u_2}{1+u_2}$ determine the inverse
map.  Next, note that $\tilde g=g+g\odot g$, where $g\odot g:S^2(W)\ra
S^2(W)$ is the map given by $(g\odot g)(ab)=g(a)g(b)$.  In our case, we
compute
\begin{equation*}
(g\odot g)(f_1f_2)=g(f_1)g(f_2)=(f_1+u_1f_2)(f_2+u_2f_2)=(1+u_2)(f_1f_2),
\end{equation*}
so that $g\odot g=(1+u_2)I\odot I$. Let us denote $\lambda(\th11)$
by $\tth11$ and $\lambda(\th12)$ by $\tth12$, so that
$\lambda_*(d_1)=\ps1+\pha11\tth11+\pha21\tth12$. Then
\begin{equation*}
\lambda_*(d_1)\tilde g=
(1+u_2)\ps1+\pha11(\tth11+u_1\tth12)+\pha21(1+u_2)\tth12.
\end{equation*}
Finally, we compute that
\begin{align*}
g\inv \lambda_*(d_1)\tilde g=&
(1+u_2)\ps1+\pha11(\tth11+u_1\tth12)+\pha21(1+u_2)\tth12\\
&+v_1(1+u_2)\ps2+
\pha12v_1(\tth11+u_1\tth12)+\pha22v_1(1+u_2)\tth12.
\end{align*}
To get the correct image we must have $\th11=\tth11+u_1\tth12$,
$\th12=(1+u_2)\tth12$, $u_2=t^1$ and $u_1=-t^2$ (which finally
justifies our tacit assumption that $1+u_2$ is invertible in \A').
Thus the map $\lambda$ is given by
$$\lambda(\th12)=\frac1{1+t^1}\th12, \qquad
\lambda(\th11)=\th11+\frac{t^2}{1+t^1}\th12.$$  Note that this map is
not even polynomial!

We computed the morphism $\lambda$ from $\A$ to
the formal algebra $\A'$, but the same ideas can be used to compute
the morphism to the algebra associated with the $n$-th order deformation.
In particular, one can compute the infinitesimal morphism. But all of
these morphisms are just the mod $\m^k$ reductions of the morphism
$\lambda$ given above.

Even though our claim that any nonzero linear combination of the cochains
$\ps1$ and $\ps2$ is equivalent as an \linf\ algebra structure to $d=\ps1$,
is just a classical fact about Lie algebras, it can be demonstrated easily
using our methodology involving a coalgebra morphism $g=I+u_1\pha12+u_2\pha22$.
If we compute $g\inv\circ d\circ\tilde g$ as before, we obtain
simply $g^*(d)=(1+u_2)\ps1-u_2\ps2$, which means that any combination
for which the coefficient of $\ps1$ is nonzero can be obtained (recall that
$1+u_2$ needs to be invertible). It is also easy to see that if
we take $g=\pha 12+\pha21$, then $g\inv=g$, $\tilde g=g+I\odot I$,
and $g^*(\ps1)=\ps2$.

\section{Deformations of a $1|1$ dimensional space}
If $W=\langle e,f\rangle$ is a $1|1$ dimensional space,
then $S^k(W)=\langle e^k,e^{k-1}f\rangle$ is a $1|1$ dimensional
space as well, and $L_k=\langle \pha ke, \pha kf, \psa ke \psa kf,\rangle$
is a $2|2$ dimensional space, where
\begin{align*}
&\pha ke(e^k)=k!e\quad&
\psa ke(e^{k-1}f)=(k-1)!e\\
&\pha kf(e^{k-1}f)=(k-1)!f\quad&
\psa kf(e^k)=k!f.\\
\end{align*}
The bracket structure is given by the table below.
$$
\begin{array}{lll}
\br{\pha me}{\pha ne}=(m-n)\pha{m+n-1}e&
\br{\pha me}{\pha nf}=(1-n)\pha{m+n-1}f\\
\br{\pha mf}{\pha nf}=0\\
\br{\pha me}{\psa ne}=(m-n+1)\psa{m+n-1}e&
\br{\pha mf}{\psa ne}=-\psa{m+n-1}e\\
\br{\pha me}{\psa nf}=-n\psa{m+n-1}f&
\br{\pha mf}{\psa nf}=\psa{m+n-1}f\\
\br{\psa me}{\psa nf}=\pha{m+n-1}e+n\pha{m+n-1}f\\
\br{\psa ne}{\psa ne}=0&
\br{\psa nf}{\psa nf}=0\\
\end{array}
$$
Let us first take care of the simple case when $d=0$.
Since every cochain is a cocycle and there are no coboundaries,
the universal infinitesimal extension is given by
\begin{equation}
d_1=\pha ne\th en +\pha nf \th fn +\psa ne \ts en +\psa nf \ts fn.
\end{equation}
This is also the formula for a miniversal deformation, but then
the following relations on the base are necessary.
\begin{align*}
&\sum_{k+l=n+1}(k-l)\th ek\th el +\ts ek\ts fl=0&\quad
&\sum_{k+l=n+1}-(k-l+1)\th ek\ts el +\th fk\ts el=0\\
&\sum_{k+1=n+1}-(l-1)\th ek\th fl +l\ts ek\ts fl=0&\quad
&\sum_{k+l=n+1}\th ek\ts fl-\th fk\ts fl=0
\end{align*}

Having disposed of the trivial codifferential, we now consider the
general case $d=r_n\psa ne +s_n\psa nf$ for some constants
$r_n,s_n\in\k$. Then a computation of $[d,d]$ yields
\begin{equation*}
\frac12\br dd=
\sum_{n=1}^\infty
\sum_{k+l=n+1}
r_ks_l(\pha ne +l\pha nf),
\end{equation*}
from which it follows that $\sum_{k+l=n+1}r_ks_l=0$ for all $n$.
Let $M$ be the least integer such that not both $r_M$ and $s_M$
vanish. If $r_M\ne 0$, then it is easy to see that this forces
$s_l=0$ for all $l$. A similar result holds if $s_M\ne 0$. Thus
the only possibilities for nontrivial codifferentials are of the
form $d=r_k\psa ke$ and $d=s_k\psa kf$.

Let us consider the case $d=r_k\psa ke$ first. Let $L$ be the
least coefficient such that $r_L\ne 0$, and for simplicity, assume
that $r_L=1$. Let $d^L=\psa kL$. Note that $\psa
{L+k}e=[d^L,\pha{k+1}f]$, so that $d=d^L+[d^L,\sum_{k=1}^\infty
r_{k+L}\psa {k+1}f]$. This resembles the case of an infinitesimal
deformation, suggesting that $d$ is equivalent to $d^L$.  Let
$\tph kf$ represent the coderivation associated to $\pha kf$
(recall that $\pha kf\in\hom(S(W),W)$).  Then the bracket is given
by $[d^L,\pha kf]=d^L\circ\tph kf$. This is because $d^L$ vanishes
on terms that do not contain an $f$, but outputs no $f$, while
$\pha kf$ requires an input of an $f$.  So in the decomposition of
the bracket, only one of the terms survives.

Let $\tl=\sum_{k=1}^\infty r_{k+L}\tph{k+1}f$. Then $\tl$ is a
coderivation of $S(W)$, so that
$$
\Delta\circ \tl=(\tl\tns I+I\tns \tl)\circ \Delta.
$$
Let $\tilde g=I+\tl$. Then we claim  that $\tilde g$ is an automorphism
of $S(W)$. For that to be true, the identity that needs to be satisfied
is
$$
\Delta\circ \tl=(\tl\tns I+I\tns \tl+\tl\tns\tl)\circ \Delta,
$$
but this holds since $(\tl\tns\tl)\circ\Delta=0$.  To see this
fact, note that $\tl$ vanishes on an element in $S(W)$ not
containing an $f$, while only one of the two terms output by
$\Delta$ could contain an $f$, since $f^2=0$. If $\tl$ was an
infinitesimal derivation, then ${\tilde g}\inv$ would simply be
$I-\tl$. This is not true here.  Nevertheless, we can compute
${\tilde g}\inv$ exactly. For simplicity, first consider the case
$L=1$ and $\tl=\tph2f$. Now,
$$
\tl\pha kf(e^nf)=(k-1)!{n\choose{k-1}}e^{n-k}f,
$$
so $\tilde g(e^n f)=e^nf +{n\choose 1}e^{n-1}f$.
It follows that
\begin{multline*}
\pha kf(e^nf+{n\choose 1}e^{n-1})=\\
(k-1)!{n\choose{k-1}}e^{n-k}f+
(k-1)!{n-1\choose{k-1}}{n\choose 1}e^{n-k-1}f.
\end{multline*}
Since $k!{n\choose k}=(k-1)!{n-1\choose k-1}{n\choose 1}$,
it is easily seen that ${\tilde g}\inv=I-\sum_{k=2}^\infty\s i\pha kf$.
For arbitrary $L$ and $\tl$, we can express
${\tilde g}\inv=I+\sum_{k=2}^\infty c_k\pha kf$, for some constants
$c_k$.

 From the consideration given above on how $\pha kf$ interacts with
$d^L$, we see that
$d={\tilde g}\inv\circ d^L\circ \tilde g$; in other words, $d$ and $d^L$
determine equivalent  \linf\ structures.

A similar argument applies to the case of $s_k\psa kf$, except
that in this case, one sees that $\psa kf\circ\tph lf=0$, because
$\pha lf$ always inputs an $f$ and outputs one, while $\psa kf$
inputs no $f$ and outputs one, so the composition is zero.  Thus
the bracket should be computed in the form $d=\tilde g\circ \psa
Lf\circ {\tilde g}\inv$ for an appropriate choice of $\tilde g$.

The \linf\ structures determined by $\psa ke$ and $\psa lf$ can never
be equivalent. To see this, note that if we extend these maps to
coderivations, then we have
$$
\begin{array}{ll}
\tps kf(e^n)={n\choose k}e^{n-k}f&
\tps kf(e^nf)=0\\
\tps ke(e^nf)={n\choose k-1}e^{n-k+2}&
\tps ke(e^n)=0\\
\end{array}.
$$
Suppose that $\tps lf={\tilde g}\inv\circ\tps ke\circ\tilde g$.
Then $\tilde g\circ\tps lf=\tps ke\circ\tilde g$. Since $\tilde
g\circ\tps lf(e^nf)=0$, it follows that $\tilde g(e^nf)$ must lie
in the kernel of $\tps ke$. But this implies that $\tilde g(e^n
f)$ would be even. This is impossible since an automorphism is an
even map, and $e^nf$ is odd. Thus no such equivalence is possible.

Next, we will compute the cohomology
$H(d_k)$ for $d_k=\pha ke$, and show that it has
dimension $2L-2$.  From this, it follows immediately that if $k\ne l$,
then $d_k$ is not equivalent to $d_l$.  A similar result holds for
$d_k=\pha kf$.  Thus there are two infinite families of nonequivalent
\linf\ structures on $W$. It is very interesting that such a simple
space gives rise to so many distinct \linf\ structures.
\subsection{The case $d^L=\psa Le$}

First, note that $D_L(\psa ke)=0$, so $\psa ke$ is a cocycle for all
$k$. No element of order less than $L$ can be a coboundary, because
the bracket of anything with $d^L$ has order at least $L$. Since
$D_L(\pha kf)=\psa {m+L-1}f$, $\psa ke$ is a coboundary if $k\ge L$.
Since $\psa kf$ is never a cocycle, it is clear that the odd cocycles
are spanned by $\psa ke$, and the dimension of the odd part of the
cohomology of $d^L$ is $L-1$.

Let $h^k=\pha ke +(k-L+1)\pha kf$. Then $D_L(h^k)=0$, and moreover
$D_L(\psa{k-L+1}f)=h^k$, when $k\ge L$, so $h^k$ is a coboundary
precisely when $k\ge L$.  Assume that $\ph=r_k\pha ke+s_k\pha kf$ is
an arbitrary even cocycle. Then
$\ph-r_kh^k=(s_k-r_k(k-L+1))\pha kf$, so because $D_L(\pha kf)\ne 0$,
it follows that $s_k=r_k(k-L+1)$ for all $k$, and thus $\ph=r_kh^k$.
The cocycles $h^k$ form a basis of the even cocycles, so
the dimension of the even part of the cohomology of $d^L$ is also
$L-1$.  Thus $\dim(H(d^L))=2L-2$.

If we consider $d^L=\psa Lf$ instead of $\psa Le$, then the odd cocycles
would be given by $\psa kf$, with those where $k\ge L$ being coboundaries,
while the even cocycles would be given by $h^k=\pha ke+L\pha kf$,
$D_L(\psa {k-L+1}e)=h^k$ for $k\ge L$. Thus we obtain the same dimension
$2L-2$ for the cohomology determined by $\psa Lf$.

Now let us work with the case $d^L=\psa Le$ and calculate a miniversal
deformation for some small values of $L$.
It will prove useful to have a table of some brackets of $h^k$ with
certain cochains.
$$
\begin{array}{ll}
\br{h^k}{h^l}=(k-l)h^{k+l-1}&
\br{h^k}{\psa le}=(L-l)\psa {k+l-1}e\\
\br{h^k}{\pha lf}=(1-l)\pha{k+l-1}f&
\br{h^k}{\psa lf}=(k-l-L+1)\psa{k+l-1}f
\end{array}
$$
\subsection{The case $d=\psa 1e$}
If $L=1$, then since $H(d^1)=0$, the miniversal deformation is just
$d^1$. This corresponds to the fact that $W$ has a differential
equipping it with the structure of a differential graded vector space,
and this differential is essentially unique.

\subsection{The case $d=\psa 2e$}
If $L=2$, then the universal infinitesimal deformation of $d^1$ is
$d^1_1=\psa 2e+h^1\theta_1+\psa 1et_1$, and
$\br{d^1_1}{d^1_1}=-\psa 1e\theta_1t_1$. Thus the infinitesimal deformation
is miniversal, with the relation $\theta_1t_1=0$. Note that in this
case, $\psa 2e$ determines a nontrivial \zt-graded Lie algebra structure on the
parity reversion of $W$, and $\psa 1e$ corresponds to
the fact that the differential on $W$
is a derivation of this Lie algebra structure.
Thus the
miniversal deformation gives a deformation of the Lie algebra structure
into an \linf\ algebra by recording the graded derivations of the Lie
algebra structure.  The part contributed by $h^1$ corresponds to a
\zt-graded antiderivation of the Lie algebra, which becomes a true
derivation only when it is multiplied by an odd parameter.

\subsection{The case $d=\psa 3e$}
For $L=3$ the universal infinitesimal deformation of $d^3$ is
\begin{equation*}
d^3_1=\psa 3e +h^1\theta_1 +h^2\theta_2+\psa 1et_1 +\psa 2et_2.
\end{equation*}
 From
\begin{equation*}
\frac12\br{d^3_1}{d^3_1}=
-h^2\theta_1\theta_2-2\psa 1e\theta_1t_1-\psa 2e(\theta_1t_2+2\theta_2t_1)
-\psa 3e\theta_2t_2
\end{equation*}
we
obtain
the
mod $\m^3$
relations
$R-2=\{\theta_1\theta_2,\theta_1t_1,\theta_1t_2+2\theta_2t_1\}$
corresponding to the coefficients of the $h^2$, $\psa 1e$ and $\psa 2e$
terms.  Since $\psa 3e$ is a coboundary, it does not give a relation,
but means that we need  to add a second order term to the codifferential.
In this example the miniversal deformation is
not given by the same formula as the universal infinitesimal deformation.
The second order deformation $d^3_2$ is given by the formula
\begin{equation*}
d^3_2=\psa 3e +h^1\theta_1 +h^2\theta_2+\psa 1et_1 +\psa 2et_2
+\pha 1f\theta_2t_2.
\end{equation*}
Finally, we compute
\begin{equation*}
\frac12\br{d^3_2}{d^3_2}=-h^2\theta_1\theta_2
-\psa 1e(2\theta_1t_1-\theta_2t_1t_2)
-\psa 2e(\theta_1t_2+2\theta_2t_1-\theta_2t^2_2).
\end{equation*}
This time no coboundary terms arise in the bracket, so that the miniversal
deformation is given by $d^3_2$, subject to the relations
$$R_3=\{\theta_1\theta_2,2\theta_1t_1-\theta_2t_1t_2,
\theta_1t_2+2\theta_2t_1-\theta_2t^2_2\}.$$
Notice that the relations come from the mod $\m^4$ relations, but
that they \emph{are modified by picking up third order terms}.
This point was misstated in \cite{ff2}, where it was declared that
any relations discovered at $n$-th order remain relations at all
orders.  What really happens is that if you reduce the
relations in the formal algebra mod $\m^{n+1}$, then they become the
$n$-th order relations.

\subsection{The case $d=\psa 4e$}
For $L=4$ the universal infinitesimal deformation of $d^4$ is
\begin{equation*}
d^4_1=\psa 4e +h^1\theta_1 +h^2\theta_2+h^3\theta^3
+\psa 1et_1 +\psa 2et_2+\psa 3et_3.
\end{equation*}
Then
\begin{align*}
\frac12\br{d^4_1}{d^4_1}=&\\
&-h^2\theta_1\theta_2
-2h^3\theta_1\theta_3
-h^4\theta_2\theta_3\\
&-3\psa 1e\theta_1t_1
-\psa 2e(2\theta_1t_2+3\theta_2t_1)
-\psa 3e(\theta_1t_3+2\theta_2t_2+3\theta_3t_1)\\
&-\psa 4e(\theta_2t_3+2\theta_3t_2)
-\psa 5e\theta_3t_3.
\end{align*}
The coefficients of $h^k$ and $\psa ke$ for $k=1,2,3$ yield 5 mod $\m^3$
relations:
\begin{equation*}
R_2=\{
\theta_1\theta_2,
\theta_1\theta_3,
\theta_1t_1,
2\theta_1t_2+3\theta_2t_1,
\theta_1t_3+2\theta_2t_2+3\theta_3t_1.
\}
\end{equation*}
There are also 3 coboundary terms, so that after adding terms to
cancel them we arrive at the following formula for the second order
deformation.
$$
\d^4_2=d^4_1+\psa 1f\theta_2\theta_3+\pha1f(\theta_2t_3+2\theta_3t_2)
+\pha2f\theta_3t_3.
$$
The bracket relation
$
\br{\psa ke}{\psa lf}=h^{k+l-1}+(L-k)\pha{k+l-1}f
$, which is easy to verify, will be used in the calculations below.
\begin{align*}
\frac12\br{d^4_2}{d^4_2}=
&h^1\theta_2\theta_3t_1
+h^2(-\theta_1\theta_2+\theta_2\theta_3t_2)
+h^3(-2\theta_1\theta_3+\theta_2\theta_3t_3)\\
&+\psa1e(-3\theta_1t_1+\theta_2t_1t_3+2\theta_3t_1t_2)\\
&+\psa2e(-2\theta_1t_2-3\theta_2t_1+\theta_2t_2t_3+2\theta_3t_2^2+\theta_3t_1t_3)
\\
&+\psa3e(-\theta_1t_3-2\theta_2t_2-3\theta_3t_1
+\theta_2t_3^2+3\theta_3t_2t_3)\\
&+3\pha1f\theta_2\theta_3t_1
+\pha2f(2\theta_2\theta_3t_2-\theta_1\theta_3t_3)\\
&+\psa4e\theta_3t_3^2
\end{align*}
The mod $\m^4$  relations are the coefficients of the $h^k$ and
$\psa ke$ terms for $k=1,2,3$, so there are six of them.  But what
about the terms that involve $\pha 1f$ and $\pha2f$?  These
cochains are not even cocycles, so what are they doing in the
expression at all? But if you check carefully, the mod $\m^4$
relations show that the coefficients are zero (mod $\m^4$). Thus
these terms do not contribute to the bracket at this level.
Finally, there is one coboundary term left, so it is necessary to
add another term, and we have
\begin{equation*}
d^4_3=d^4_2-\pha 1f\theta_3t_3^2.
\end{equation*}
After computing the boundary, the relations given by the coefficients
of $\psa 1e$, $\psa2e$ and $\psa3e$ are modified by adding the fourth
order term $-\theta_3t_it_3^2$ to the coefficient of $\psa ie$.
The coefficients of the cocycles which are not coboundaries
determine six relations on the base
which is of the form $\A=\k[[\theta_1,\theta_2,\theta_3,t_1,t_2,t_3]]/R$.
The coefficient of $\pha1f$ is one of the relations, so it is clearly
zero. It is also easy to show that the coefficient of $\pha2f$ is equal
to zero mod $\m^4$, but since we are claiming that the deformation is
miniversal, the coefficient must be exactly equal to zero.
Multiplying the coefficient of $\psa3f$ by $\theta_3$ yields
\begin{multline*}
\theta_1\theta_3t_3+2\theta_2\theta3t_2-\theta_2\theta_3t_3^2
\\=\theta_1\theta_3t_3 +2\theta_2\theta_3t_2-2\theta_1\theta_3t_3\\=
2\theta_2\theta_3t_2-\theta_1\theta_3t_3,
\end{multline*}
where we used the relation $\theta_2\theta_3t_3=2\theta_1\theta_3$
in the second step.  This is precisely the coefficient of $\pha2f$,
so the coefficient of this term is zero. (Of course, it is really not
necessary to show that these coefficients are zero, because it follows
from the general theory that this must be the case.)
Thus $d^4$ gives a miniversal deformation.
Note that all terms in the relations have odd parameters in them.
These computations illustrate the importance of introducing odd
parameters, because otherwise $d_1$ would be a miniversal
deformation.

It is not obvious how to write down the formula for the miniversal deformation
of $d^L$
in general, but there are some things which can be easily shown in the
general case.  First, there will be $2L-2$ relations, (when $L\ge 4$).
Secondly, by studying the examples so far, it can be seen that the
$(L-1)$-th order deformation is miniversal.  Thus only a finite number of
computations is necessary in order to compute the miniversal deformation.
\section{Conclusions}
In the examples we have studied, we have given explicit constructions
of miniversal deformations.  Previous works have only computed the
base of a miniversal deformation.  The main complication in these
constructions is the requirement of exact knowledge of the bracket
structure on the space of cochains.  Since most of the examples in
the literature where cohomology of Lie algebras is studied do not
contain this information, the computation of versal deformations
of these structures will require more information than is currently
easily available.  We have just touched on the beginning of the subject.

In a future work,
the authors plan to investigate  versal deformations of \linf\ algebras
with invariant inner products, which is related to the cyclic
cohomology of these algebras.  In the case of ordinary Lie algebras,
only reductive algebras have invariant inner products, so the
deformation theory is not so interesting.  For example, simple
Lie algebras have no deformations, even as \linf\ algebras.
For super Lie algebras, and \linf\ algebras
in general the picture is not so restricted.
\bibliographystyle{amsplain}
%\bibliography{global}

\begin{thebibliography}{10}

\bibitem{aksz}
M.~Alexandrov, M.~Kontsevich, A.~Schwarz, and O.~Zaboronsky, \emph{The geometry
  of the master equation and topological quantum field theory}, Internat.
  Journ. Modern Phys. \textbf{A12} (1997), 1405--1423.

\bibitem{bfls}
G.~Barnich, R.~Fulp, T.~Lada, and J.~Stasheff, \emph{The sh {L}ie structure of
  {P}oisson brackets in field theory}, Comm. Math. Phys. \textbf{191} (1998),
  no.~3, 585--601.

\bibitem{fi}
A.~Fialowski, \emph{Deformations of {L}ie algebras}, Mathematics of the
  USSR-Sbornik \textbf{55} (1986), no.~2, 467--473.

\bibitem{fi2}
\bysame, \emph{An example of formal deformations of {L}ie algebras}, NATO
  Conference on Deformation Theory of Algebras and Applications (Dordrecht),
  Kluwer, 1988, Proceedings of a NATO conference held in Il Ciocco, Italy,
  1986, pp.~375--401.

\bibitem{ff3}
A.~Fialowski and D.~Fuchs, \emph{Singular deformations of Lie algebras on an
  example}, Topics in Singularity Theory (Providence, RI) (A.~Varchenko and
  V.~Vassilie, eds.), A.M.S. Translation Series 2, Vol.180, Amer.\ Math.\ Soc.,
  1997, V.~I.~Arnold {\rm 60}${}^{th}$ Anniversary Collection,.

\bibitem{ff2}
\bysame, \emph{Construction of miniversal deformations of {L}ie algebras},
  Journal of Functional Analysis (1999), no.~161(1), 76--110.

\bibitem{fiapen1}
A.~Fialowski and M.~Penkava, \emph{Deformation theory of infinity algebras},
  preprint, 2000.

\bibitem{fp}
A~Fialowski and G~Post, \emph{Versal deformation of the Lie algebra
  $\mathbf{L_2}$}, Journal of Algebra (2001) no.~236(1).

\bibitem{kon}
M.~Kontsevich, \emph{Feynman diagrams and low dimensional topology}, First
  European Congress of Mathematics, Paris, 1992, Birkhauser, Basel, 1994,
  pp.~97--121.

\bibitem{lm}
T.~Lada and M.~Markl, \emph{Strongly homotopy {L}ie algebras}, Comm. in Algebra
  \textbf{23} (1995), 2147--2161.

\bibitem{ls}
T.~Lada and J.~Stasheff, \emph{Introduction to sh {L}ie algebras for
  physicists}, Intern. J. Theor. Phys \textbf{32} (1993), 1087--1103, Preprint
  hep-th 9209099.

\bibitem{mar2}
M.~Markl, \emph{Cyclic operads and the homology of graph complexes}, Preprint
  Math/9801095, 1996.

\bibitem{pen1}
M.~Penkava, \emph{\hbox{$L_\infty$} algebras and their cohomology}, Preprint
  q-alg 9512014, 1996.

\bibitem{pen3}
\bysame, \emph{Infinity algebras and the homology of graph complexes}, Preprint
  q-alg 9601018, 1996.

\bibitem{ps2}
M.~Penkava and A.~Schwarz, \emph{\hbox{$A_\infty$} algebras and the cohomology
  of moduli spaces}, Dynkin Seminar, vol. 169, American Mathematical Society,
  1995, pp.~91--107.

\bibitem{lp}
M.~Penkava and L.~Weldon, \emph{Infinity algebras, {M}assey products, and
  deformations}, Preprint math/9808058, 1996.

\bibitem{ss}
M.~Schlessinger and J.~Stasheff, \emph{The {L}ie algebra structure of tangent
  cohomology and deformation theory}, Journal of Pure and Applied Algebra
  \textbf{38} (1985), 313--322.

\bibitem{sta1}
J.D. Stasheff, \emph{On the homotopy associativity of {H}-spaces {I}},
  Transactions of the AMS \textbf{108} (1963), 275--292.

\bibitem{sta2}
\bysame, \emph{On the homotopy associativity of {H}-spaces {II}}, Transactions
  of the AMS \textbf{108} (1963), 293--312.

\bibitem{sta3}
\bysame, \emph{Closed string field theory, strong homotopy lie algebras and the
  operad actions of moduli spaces}, Conf. Proc. Lecture Notes Math. Phys., III,
  Internat. Press, Cambridge, MA, 1994, Perspectives in mathematical physics,
  pp.~265--288.

\end{thebibliography}
%\providecommand{\bysame}{\leavevmode\hbox to3em{\hrulefill}\thinspace}

\end{document}